\documentclass[12pt,a4paper]{article}
\usepackage{amsfonts}
\usepackage{amssymb}
\usepackage{mathrsfs}
\usepackage{amsmath}
\usepackage{amsthm}
\usepackage{enumerate}
\usepackage{cite}
\usepackage[all]{xy}
\usepackage{extarrows}
\usepackage{indentfirst}
\allowdisplaybreaks
\newtheorem{defn}{Definition}[section]
\newtheorem{thm}[defn]{Theorem}
\newtheorem{lem}[defn]{Lemma}
\newtheorem{prop}[defn]{Proposition}
\newtheorem{cor}[defn]{Corollary}
\newtheorem{eg}[defn]{Example}
\newtheorem{re}[defn]{Remark}

\newcommand{\bdefn}{\begin{defn}}
\newcommand{\edefn}{\end{defn}}
\newcommand{\bthm}{\begin{thm}}
\newcommand{\ethm}{\end{thm}}
\newcommand{\blem}{\begin{lem}}
\newcommand{\elem}{\end{lem}}
\newcommand{\bprop}{\begin{prop}}
\newcommand{\eprop}{\end{prop}}
\newcommand{\bcor}{\begin{cor}}
\newcommand{\ecor}{\end{cor}}
\newcommand{\beg}{\begin{eg}}
\newcommand{\eeg}{\end{eg}}
\newcommand{\bre}{\begin{re}}
\newcommand{\ere}{\end{re}}
\newcommand{\bpf}{\begin{proof}}
\newcommand{\epf}{\end{proof}}
\newcommand{\benu}{\begin{enumerate}}
\newcommand{\eenu}{\end{enumerate}}
\newcommand{\bc}{\begin{center}}
\newcommand{\ec}{\end{center}}
\newcommand{\bea}{\begin{eqnarray}}
\newcommand{\eea}{\end{eqnarray}}
\newcommand{\Bea}{\begin{eqnarray*}}
\newcommand{\Eea}{\end{eqnarray*}}
\newcommand{\beq}{\begin{equation}}
\newcommand{\eeq}{\end{equation}}
\newcommand{\Beq}{\begin{equation*}}
\newcommand{\Eeq}{\end{equation*}}
\newcommand{\bspl}{\begin{split}}
\newcommand{\espl}{\end{split}}

\newcommand{\ch}{{\rm ch}}
\newcommand{\id}{{\rm id}}

\newcommand{\End}{{\rm End}}
\newcommand{\Ker}{{\rm Ker}}
\newcommand{\K}{\mathbb{K}}

\newcommand{\N}{\mathbb{N}}
\newcommand{\Z}{\mathbb{Z}}
\newcommand{\I}{\mathcal{I}}
\newcommand{\J}{\mathcal{J}}
\newcommand{\M}{\mathfrak{M}}
\begin{document}
\title{{\bf  Some structure theories of Leibniz triple systems}}
\author{\normalsize \bf Yao Ma$^{1}$,  Liangyun Chen$^{2*}$}
\date{\centering {{\small{ {}$^1$School of Mathematical Sciences, University of Science and Technology of China,\\ Hefei 230026, CHINA
 \\{}$^2$School of Mathematics and Statistics,  Northeast Normal University,\\ Changchun 130024, CHINA}}}} \maketitle
\date{}

\begin{abstract}
In this paper, we investigate the Leibniz triple system $T$ and its
universal Leibniz envelope $U(T)$. The involutive automorphism of
$U(T)$ determining $T$ is introduced, which gives a characterization
of the $\Z_2$-grading of $U(T)$. We give the relationship between
the solvable radical $R(T)$ of $T$ and $Rad(U(T))$, the solvable
radical of $U(T)$. Further, Levi's theorem for Leibniz triple
systems is obtained. Moreover, the relationship between the nilpotent radical of $T$ and that of $U(T)$ is studied. Finally, we introduce the notion
 of representations of a Leibniz triple system.

\bigskip

\end{abstract}
\renewcommand{\thefootnote}{\fnsymbol{footnote}}
\footnote[0]{* Corresponding author(L. Chen):
chenly640@nenu.edu.cn.} \footnote[0]{Supported by  NNSF of China
(Nos. 11171055 and 11471090) and China Postdoctoral Science
Foundation(No. 159347). }

\section{Introduction}
The notion of Leibniz algebras was introduced by Loday \cite{L2}, which is a ``nonantisymmetric'' generalization of Lie algebras. A Leibniz algebra is a
 vector space equipped with a bilinear bracket satisfying the (right) Leibniz identity
\beq\label{Leibniz alg}
[[x, y], z]=[[x, z], y]+[x, [y, z]].
\eeq
A Leibniz algebra whose bracket is antisymmetric is a Lie algebra. In fact, such algebras were considered by Bloh in 1965, who called them $D$-algebras \cite{B2}. Later,
Loday introduced this class of algebras to search an ``obstruction'' to the periodicity in algebraic $K$-theory. So far many results of Leibniz algebras have been studied,
including extending important theorems in Lie algebras: there are analogs of Lie's theorem, Engel's theorem, Levi's theorem and Cartan's criterion for Leibniz
algebras \cite{AO,AAO,AAO2,B3,B4,DMS,GKO,G,LP,P}.

Leibniz triple systems were introduced by Bremner and S\'{a}nchez-Ortega \cite{BS}. Leibniz triple systems are defined in a functorial manner using
the Kolesnikov-Pozhidaev algorithm, which takes the defining identities for a variety of algebras and produces the defining identities for the corresponding
 variety of dialgebras \cite{K}. As an application of this algorithm, it is showed that associative dialgebras and Leibniz algebras can be obtained from associative and Lie algebras, respectively. In \cite{BS}, Leibniz triple systems were obtained by applying the Kolesnikov-Pozhidaev algorithm to Lie triple systems. Furthermore, Leibniz triple systems are related to Leibniz algebras in the same way that Lie triple systems related to Lie algebras. So it is natural to prove analogs of results from the theory of Lie triple systems to Leibniz triple systems.

This paper proceeds as follows. Section 2 devotes to some basic
facts about a Leibniz triple system $T$. In Section 3, we give the
involutive automorphism of $U(T)$ determining $T$, where $U(T)$ is
the universal Leibniz envelope of $T$, using which to investigate
the connection between automorphisms of $T$ and those of $U(T)$, and
to describe the $\Z_2$-graded subspace of $U(T)$. In Section 4, the
solvable radical of a Leibniz triple system and the definition of a
semisimple Leibniz triple system are introduced. In Section 5,
Levi's theorem is extended to the case of Leibniz triple systems. In Section 6, the notion of nilpotent ideals of $T$ is introduced. Section 7
is devoted to the representations of a Leibniz triple system.

In the sequel, all vector spaces are finite-dimensional and defined
over a fixed but arbitrarily chosen field $\K$ of characteristic
$0$.

\section{Preliminaries}
\bdefn{\rm\cite{L,M}}
A vector space $T$ together with a trilinear map $(a, b, c)\mapsto[abc]$ is called a \textbf{Lie triple system} if
\begin{gather*}
[aac]=0,\\
[abc]+[bca]+[cab]=0,\\
[ab[cde]]=[[abc]de]+[c[abd]e]+[cd[abe]],
\end{gather*}
for all $a,b,c,d,e\in T$.
\edefn

\bdefn{{\rm\cite{BS}}}
A \textbf{Leibniz triple system} (LeibTS for short) is a vector space $T$ with a trilinear operation $\{\cdot,\cdot,\cdot\}: T\times T\times T\rightarrow T$ satisfying
\begin{gather}
\{a\{bcd\}e\}=\{\{abc\}de\}-\{\{acb\}de\}-\{\{adb\}ce\}+\{\{adc\}be\},\label{VIP1}\\
\{ab\{cde\}\}=\{\{abc\}de\}-\{\{abd\}ce\}-\{\{abe\}cd\}+\{\{abe\}dc\}.\label{VIP2}
\end{gather}
\edefn

A Lie triple system gives a LeibTS with the same ternary product.
If $L$ is a Leibniz algebra with product $[\cdot,\cdot]$, then $L$ becomes a LeibTS by putting $\{xyz\}=[[xy]z]$.
The following examples of LeibTSs are induced by Leibniz algebras.
\beg
(1) Let $T=\K\langle x, y\rangle$ with multiplications given by $\{xyy\}=\{yyy\}=x$, and the rest are 0. Then $T$ is a LeibTS.

(2) Let $T=\K\langle x, y, z\rangle$ with multiplications given by $\{xzz\}=\{yzz\}=\{zzz\}=x$, and the rest are 0. Then $T$ is a LeibTS.

(3) Let $T=\K\langle x, a, b, c, d\rangle$ with the following multiplications
      \begin{gather*}
      \{axx\}=-\{xax\}=a, \quad \{bax\}=\{cxx\}=c,\\
      \{abx\}=\{xab\}=\{axb\}=\{xcx\}=\{xdx\}=\{dxx\}=d,
      \end{gather*}
      and the rest are 0. Then $T$ is a LeibTS.

More examples refer to \cite{BS}.
\eeg

\bprop{{\rm\cite{BS}}}
Let $T$ be a LeibTS. Then the following identities hold.
\begin{gather}
\{a\{bcd\}e\}+\{a\{cbd\}e\}=0,~~ \{ab\{cde\}\}+\{ab\{dce\}\}=0,\label{comm}\\
\{a\{bcd\}e\}\!+\{a\{cdb\}e\}\!+\{a\{dbc\}e\}\!=0,~~
\{ab\{cde\}\}\!+\{ab\{dec\}\}\!+\{ab\{ecd\}\}\!=0,\label{cyclic}\\
\{\{cde\}ba\}-\{\{cde\}ab\}-\{\{cba\}de\}+\{\{cab\}de\}-\{c\{abd\}e\}-\{cd\{abe\}\}=0.\label{VIP3}
\end{gather}
\eprop

\bre
Any two of Equations (\ref{VIP1}), (\ref{VIP2}) and (\ref{VIP3}) induce the rest one.
\ere

\bthm{{\rm\cite{BS}}}\label{uni Leib envelop}
Suppose that $T$ is a LeibTS with product $\{xyz\}$. The \textbf{universal Leibniz envelope} of $T$ is the
quotient algebra $U(T)=A(T)/I(T)$ with $A(T)$ the free Leibniz algebra and $I(T)$ the ideal generated by the
elements $[[x, y], z]-\{xyz\}$. Then $U(T)=T\oplus(T\otimes T)$ together with product $[\cdot, \cdot]$ satisfies
\begin{gather*}
[a, b]=a\otimes b,~~ [a\otimes b, c]=\{abc\},~~ [a, b\otimes c]=\{abc\}-\{acb\},\\
[a\otimes b, c\otimes d]=\{abc\}\otimes d-\{abd\}\otimes c.
\end{gather*}
It is straightforward to verify that $U(T)$ is a Leibniz algebra and $T\otimes T$ is a subalgebra of $U(T)$.
\ethm

Generally speaking, $U(T)$ is not a Lie algebra even if $T$ is a Lie triple system, so we try to add a condition.

\bprop
Let $T$ be a Lie triple system and $U(T)$ be defined as in Theorem \ref{uni Leib envelop}. Denote by $K(T)=\{\sum a\otimes b\in T\otimes T\ |\ [\sum a\otimes b, T] =[T, \sum a\otimes b]=0\}$. If $K(T)=0$, then $U(T)$ is a Lie algebra.
\eprop
\bpf
We need only to show that for all $x=a+\sum b\otimes c\in U(T)$, $[x, x]=0$, since $U(T)$ is already a Leibniz algebra by Theorem \ref{uni Leib envelop}.

Note that $T$ is a Lie triple system. Take $a, b, c\in T$, then one has
$$[a\otimes a, b]=\{aab\}=0,\quad [b, a\otimes a]=\{baa\}-\{baa\}=0,$$
which imply $[a,a]=a\otimes a\in K(T)=0$. It follows from (\ref{VIP2}) and the definition of Lie triple systems that
\begin{gather*}
\begin{aligned}
&[[a\otimes b,a\otimes b], c]=[\{aba\}\otimes b-\{abb\}\otimes a, c] =\{\{aba\}bc\}-\{\{abb\}ac\}\\
=&\{ab\{abc\}\}+\{\{abc\}ab\}-\{\{abc\}ba\} =\{ab\{abc\}\}+\{\{abc\}ab\}+\{b\{abc\}a\}=0,
\end{aligned}
\intertext{and}
\begin{aligned}
&[c, [a\otimes b,a\otimes b]]=[c, \{aba\}\otimes b-\{abb\}\otimes a]\\
=&\{c\{aba\}b\}-\{cb\{aba\}\}-\{c\{abb\}a\}+\{ca\{abb\}\} \\ =&\{b\{aba\}c\}+\{\{abb\}ac\}=\{ab\{bac\}\}-\{ba\{abc\}\}=0,
\end{aligned}
\end{gather*}
that is, $[a\otimes b,a\otimes b]\in K(T)=0$.

Now take $a+\sum_{i=1}^{n}b_i\otimes c_i\in U(T)$, then
\begin{align*}
&\left[a+\sum_{i=1}^{n}b_i\otimes c_i, a+\sum_{i=1}^{n}b_i\otimes c_i\right]\\
=&a\otimes a+\sum_{i=1}^n(\{b_ic_ia\}+\{ab_ic_i\}-\{ac_ib_i\}) +\sum_{i,j=1}^n[b_i\otimes c_i, b_j\otimes c_j]\\
=&\sum_{i<j}([b_i\otimes c_i, b_j\otimes c_j]+[b_j\otimes c_j, b_i\otimes c_i]).
\end{align*}
Note that for all $d\in T$, $[d, [b_i\otimes c_i, b_j\otimes c_j]+[b_j\otimes c_j, b_i\otimes c_i]]=0$ by using (\ref{Leibniz alg}), and
\begin{align*}
&[[b_i\otimes c_i, b_j\otimes c_j], d]=[\{b_ic_ib_j\}\otimes c_j -\{b_ic_ic_j\}\otimes b_j, d]\\
=&\{\{b_ic_ib_j\}c_jd\}-\{\{b_ic_ic_j\}b_jd\} =\{b_ic_i\{b_jc_jd\}\}-\{b_jc_j\{b_ic_id\}\}.
\end{align*}
Then $[[b_i\otimes c_i, b_j\otimes c_j]+[b_j\otimes c_j, b_i\otimes c_i], d]=0$,
and so
$$[b_i\otimes c_i, b_j\otimes c_j]+[b_j\otimes c_j, b_i\otimes c_i] \in K(T)=0.$$ Thus $[x, x]=0$.
\epf

\bdefn
Let $I$ be a subspace of a LeibTS $T$. Then $I$ is called a \textbf{subsystem} of $T$, if $\{III\}\subseteq I$; $I$ is called an \textbf{ideal} of $T$,
if $\{ITT\}+\{TIT\}+\{TTI\}\subseteq I$.
\edefn

Let $I$ be an ideal of a LeibTS $T$ and define some notations as follows.
\begin{gather*}
i(I)=I\otimes T+T\otimes I,\\
j(I)=\left\{\sum a\otimes b\in T\otimes T\ \Big{|}\ \left[\sum a\otimes b, T\right]+\left[T, \sum a\otimes b\right]\subseteq I\right\},\\
\I(I)=I+i(I),\quad \J(I)=I+j(I).
\end{gather*}
It is easy to see that $K(T)=j(0)=\J(0)$. The notion of $K(T)$ will be used in Section 4, too.
\bthm
If $I$ is an ideal of a LeibTS $T$, then
\benu[(1)]
\item $i(I)\subseteq j(I)$,  $\I(I)\subseteq\J(I)$.
\item $i(I)$ and $j(I)$ are ideals of $T\otimes T$.
\item $\I(I)$ and $\J(I)$ are ideals of $U(T)$.
\eenu
\ethm
\bpf
(1) follows from
\begin{gather*}
[i(I), T]=[I\otimes T+T\otimes I, T]\subseteq \{ITT\}+\{TIT\}\subseteq I
\intertext{and}
[T, i(I)]=[T, I\otimes T+T\otimes I]\subseteq\{TIT\}+\{TTI\}\subseteq I.
\end{gather*}

(2) $i(I)$ is an ideal of $T\otimes T$ since
\begin{gather*}
[i(I), T\otimes T]=[I\otimes T+T\otimes I, T\otimes T]
\subseteq\{ITT\}\otimes T+\{TIT\}\otimes T\subseteq I\otimes T\subseteq i(I)
\intertext{and}
[T\otimes T, i(I)]=[T\otimes T, I\otimes T\!+\!T\otimes I]\subseteq\!\{TTI\}\otimes T\!+\!\{TTT\}\otimes I\subseteq\! I\otimes T\!+\!T\otimes I\!=i(I).
\end{gather*}

To show that $j(I)$ is an ideal of $T\otimes T$, we take $\sum t\otimes s\in T\otimes T$ and $\sum a\otimes b\in j(I)$. It is sufficient to prove that for all $c\in T$,
\begin{gather*}
\left[\left[\sum a\otimes b, \sum t\otimes s\right], c\right]\in I,\quad \left[c, \left[\sum a\otimes b, \sum t\otimes s\right]\right]\in I,\\
\left[\left[\sum t\otimes s, \sum a\otimes b\right], c\right]\in I,\quad \left[c, \left[\sum t\otimes s, \sum a\otimes b\right]\right]\in I.
\end{gather*}

In fact, it follows from (\ref{VIP2}) that
\begin{gather*}
\begin{aligned}
&\left[\left[\sum a\otimes b, \sum t\otimes s\right], c\right] =\sum(\left[\{abt\}\otimes s-\{abs\}\otimes t, c\right])\\
=&\sum(\{\{abt\}sc\}-\{\{abs\}tc\})
=\sum(\{[a\otimes b, t]sc\}-\{[a\otimes b, s]tc\})\in I,
\end{aligned}\\
\begin{aligned}
&\left[c, \left[\sum a\otimes b, \sum t\otimes s\right]\right] =\sum[c, \{abt\}\otimes s-\{abs\}\otimes t]\\
=&\sum(\{c\{abt\}s\}-\{cs\{abt\}\}-\{c\{abs\}t\}+\{ct\{abs\}\})\\
=&\sum(\{c[a\otimes b, t]s\}-\{cs[a\otimes b, t]\}-\{c[a\otimes b, s]t\}+\{ct[a\otimes b, s]\})\in I,
\end{aligned}
\intertext{and}
\begin{aligned}
&\left[\left[\sum t\otimes s, \sum a\otimes b\right], c\right] =\sum(\left[\{tsa\}\otimes b-\{tsb\}\otimes a, c\right])\\
=&\sum(\{\{tsa\}bc\}-\{\{tsb\}ac\})
=\sum(\{ts\{abc\}\}+\{\{tsc\}ab\}-\{\{tsc\}ba\})\\
=&\sum(\{ts[a\otimes b, c]\}+[\{tsc\}, a\otimes b])\in I.
\end{aligned}
\end{gather*}
Moreover, by (\ref{Leibniz alg}), $\left[c, \left[\sum t\otimes s, \sum a\otimes b\right]\right] =-\left[c, \left[\sum a\otimes b, \sum t\otimes s\right]\right]\in I$.

(3) It is routine to prove that $\I(I)$ is an ideal of $U(T)$ by the facts of $I$ being an ideal of $T$ and $i(I)$ being
 an ideal of $T\otimes T$. $\J(I)$ is an ideal of $U(T)$ since
\begin{align*}
[\J(I), U(T)]=&[I+j(I), T+T\otimes T]
=I\otimes T+[j(I), T]+[I, T\otimes T]+[j(I), T\otimes T]\\
\subseteq&i(I)+I+j(I)=I+j(I)=\J(I),
\intertext{and}
[U(T), \J(I)]=&[T+T\otimes T, I+j(I)]
=T\otimes I+[T\otimes T, I]+[T, j(I)]+[T\otimes T, j(I)]\\
\subseteq&i(I)+I+j(I)=\J(I).
\end{align*}
Thus the proof is completed.
\epf

\section{The involutive automorphism of $U(T)$ determining $T$}
Theorem \ref{uni Leib envelop} says that $U(T)$ is a $\Z_2$-graded Leibniz algebra with $U(T)_{\bar{0}}=T\otimes T$ and $U(T)_{\bar{1}}=T$. Then every Leibniz triple
system $T$ is the ${\bar{1}}$ component of a $\Z_2$-graded Leibniz algebra. Moreover, for a $\Z_2$-graded Leibniz algebra $L=L_{\bar{0}}\oplus L_{\bar{1}}$, $L_{\bar{1}}$
has the structure of a Leibniz triple system.

The $\Z_2$-grading of $U(T)$ can be characterized by an involutive automorphism. In what follows we will use this involutive automorphism of $U(T)$ to investigate the
connection between automorphisms of $T$ and those of $U(T)$, and describe the $\Z_2$-graded subspace of $U(T)$.
\bthm\label{theta}
Let $T$ be a LeibTS and $U(T)$ the universal Leibniz envelope of $T$. Then there is a unique involutive automorphism $\theta: U(T)\rightarrow U(T)$ such that $$T=\{x\in U(T)\ |\ \theta(x)=-x\}.$$
\ethm
\bpf
For $a+\sum b\otimes c\in U(T)$, define $\theta\left(a+\sum b\otimes c\right)=-a+\sum b\otimes c$. It is clear that $\theta$ is linear and $\theta^2=\id$. Note that
\begin{align*}
&\left[\theta\left(a+\sum b\otimes c\right), \theta\left(a'+\sum b'\otimes c'\right)\right]=\left[-a+\sum b\otimes c, -a'+\sum b'\otimes c'\right]\\
=&a\otimes a'-\sum\{bca'\}-\sum\{ab'c'\}+\sum\{ac'b'\}+\sum\{bcb'\}\otimes c'-\sum\{bcc'\}\otimes b',
\intertext{and}
&\theta\left(\left[a+\sum b\otimes c, a'+\sum b'\otimes c'\right]\right)\\
=&\theta\left(a\otimes a'+\sum\{bca'\}+\sum\{ab'c'\}-\sum\{ac'b'\}+\sum\{bcb'\}\otimes c'-\sum\{bcc'\}\otimes b'\right)\\
=&a\otimes a'-\sum\{bca'\}-\sum\{ab'c'\}+\sum\{ac'b'\}+\sum\{bcb'\}\otimes c'-\sum\{bcc'\}\otimes b'.
\end{align*}
Then $\theta$ is an involutive automorphism of $U(T)$ and $T=\{x\in U(T)\ |\ \theta(x) =-x\}.$

If $\phi$ is another involutive automorphism of $U(T)$ with $T=\{x\in U(T)\ |\ \phi(x)=-x\}$, then $\theta|_T=\phi|_T=-\id_T$ and
\begin{align*}
&\theta\left(a+\sum b\otimes c\right)\\
=&\theta\left(a+\sum [b, c]\right)
=\theta(a)+\sum[\theta(b), \theta(c)]=\phi(a)+\sum[\phi(b), \phi(c)]
=\phi\left(a+\sum [b, c]\right)\\=&\phi\left(a+\sum b\otimes c\right),
\end{align*}
which implies $\theta=\phi$, so $\theta$ is unique. $\theta$ is called the \textbf{involutive automorphism of} $U(T)$ \textbf{determining} $T$.
\epf

\bprop
Suppose that $T$ is a LeibTS and $f$ is an automorphism of $T$. Then there is a unique automorphism $\tilde{f}$ of $U(T)$ such that $f=\tilde{f}|_T$ and $\tilde{f}\theta=\theta\tilde{f}$, where $\theta$ is the involutive automorphism of $U(T)$ determining $T$.
\eprop
\bpf
Define $\tilde{f}: U(T)\rightarrow U(T)$ by $\tilde{f}(a+\sum b\otimes c)=f(a)+\sum f(b)\otimes f(c)$. It is clear that $\tilde{f}$ is well defined, $\tilde{f}|_T=f$ and $\tilde{f}$ is a linear isomorphism. Moreover, we have
\begin{align*}
&\tilde{f}\left(\left[a+\sum b\otimes c, a'+\sum b'\otimes c'\right]\right)\\
=&\tilde{f}\left(a\otimes a'+\sum\{bca'\}+\sum\{ab'c'\}-\sum\{ac'b'\}+\sum\{bcb'\}\otimes c'-\sum\{bcc'\}\otimes b'\right)\\
=&f(a)\otimes f(a')+\sum\{f(b)f(c)f(a')\}+\sum\{f(a)f(b')f(c')\}-\sum\{f(a)f(c')f(b')\}\\ &+\sum\{f(b)f(c)f(b')\}\otimes f(c')-\sum\{f(b)f(c)f(c')\}\otimes f(b')\\
=&\left[f(a)+\sum f(b)\otimes f(c), f(a')+\sum f(b')\otimes f(c')\right]\\
=&\left[\tilde{f}\left(a+\sum b\otimes c\right), \tilde{f}\left(a'+\sum b'\otimes c'\right)\right],
\end{align*}
and
\begin{align*}
&\tilde{f}\theta\left(a+\sum b\otimes c\right)=\tilde{f}\left(-a+\sum b\otimes c\right)=-f(a)+\sum f(b)\otimes f(c)\\
=&\theta\left(f(a)+\sum f(b)\otimes f(c)\right)=\theta\tilde{f}\left(a+\sum b\otimes c\right).
\end{align*}
Hence $\tilde{f}$ is an automorphism and $\tilde{f}\theta=\theta\tilde{f}$. $\tilde{f}$ is unique since $U(T)$ is generated by $T$.
\epf

\bprop
If $f$ is an automorphism of $U(T)$ satisfying $f\theta=\theta f$, then $f(T)\subseteq T$ and $f|_T$ is an automorphism of $T$.
\eprop
\bpf
For all $a\in T$, $f\theta(a)=f(-a)=-f(a)$. Since $f\theta=\theta f$, it follows that $\theta f(a)=-f(a)$, which implies $f(a)\in T$ by Theorem \ref{theta}.
 It is clear that $f|_T$ is an injective endomorphism of $T$. Note that
$$f(T\otimes T)=f([T, T])=[f(T), f(T)]=f(T)\otimes f(T)\subseteq T\otimes T,$$
and $f$ is surjective, so $f|_T$ is surjective. Thus $f|_T$ is an automorphism of $T$.
\epf

Note that $U(T)$ is a $\Z_2$-graded Leibniz algebra with $U(T)_{\bar{0}}=T\otimes T$ and $U(T)_{\bar{1}}=T$. Then a subspace $V$ of $U(T)$ is $\Z_2$-graded if
and only if $V=V\cap T+V\cap(T\otimes T)$.
Suppose that $V$ is a $\theta$-invariant subspace of $U(T)$. Then for $a+\sum b\otimes c\in V$, $\theta\left(a+\sum b\otimes c\right)=-a+\sum b\otimes c \in V$.
 Since $\ch\K=0$, it follows $a\in V\cap T$ and $\sum b\otimes c\in V\cap(T\otimes T)$, that is, $V=V\cap T+V\cap(T\otimes T)$. Conversely, if $V$ is a subspace of
  $U(T)$ such that $V=V\cap T+V\cap(T\otimes T)$, then $a+\sum b\otimes c\in V$ implies that $a\in V$ and $\sum b\otimes c\in V$, so $\theta\left(a+\sum b\otimes
  c\right)=-a+\sum b\otimes c \in V$ and $V$ is a $\theta$-invariant subspace of $U(T)$.

Therefore, $V$ being a $\theta$-invariant subspace of $U(T)$ is equivalent to that $V$ is a $\Z_2$-graded subspace of $U(T)$.

\bthm\label{i(I)<I+M<j(I)}
Let $V$ be a $\theta$-invariant subspace of $U(T)$. If $V$ is an ideal of $U(T)$, then $V\cap T$ is an ideal of $T$ and $V\cap (T\otimes T)$ is an ideal
of $T\otimes T$ such that
$$i(V\cap T)\subseteq V\cap (T\otimes T)\subseteq j(V\cap T),\quad \I(V\cap T)\subseteq V\subseteq \J(V\cap T).$$
\ethm
\bpf Denote $I=V\cap T$ and $\M=V\cap (T\otimes T)$. Then $V=I+\M$. Using the condition that $V$ is an ideal of $U(T)$, one has
\begin{gather*}
[V, U(T)]=[I+\M, T+T\otimes T]=I\otimes T+[\M, T]+[I,[T, T]]+[\M, T\otimes T]\subseteq I+\M,\\
[U(T), V]=[T+T\otimes T, I+\M]=T\otimes I+[[T, T], I]+[T, \M]+[T\otimes T, \M]\subseteq I+\M.
\end{gather*}
It follows that
\begin{gather*}
I\otimes T+T\otimes I+[\M, T\otimes T]+[T\otimes T, \M]\subseteq \M,\\
[\M, T]+[I,[T, T]]+[[T, T], I]+[T, \M]\subseteq I.
\end{gather*}
Then $i(I)\subseteq\M\subseteq j(I)$ and $\M$ is an ideal of $T\otimes T$. Note that
\begin{gather*}
\{TTI\}=[[T, T], I]\subseteq I,\\
\{TIT\}=[[T, I], T]=[T\otimes I, T]\subseteq[\M, T]\subseteq I,\\
\{ITT\}=[[I, T], T]=[I\otimes T, T]\subseteq[\M, T]\subseteq I.
\end{gather*}
Thus $I$ is an ideal of $T$ and $\I(I)\subseteq I+\M\subseteq \J(I)$.
\epf

\section{The solvability of a Leibniz triple system}
Given an arbitrary Leibniz algebra $L$, the derived sequence of an ideal $I$ of $L$ is defined as
$$I^{(0)}=I, I^{(n+1)}=[I^{(n)}, I^{(n)}], \quad\forall n\geq1.$$
Suppose that $I$ is an ideal of a LeibTS $T$. Let $I^{[0]}=I$ and for $n\geq1$, $$I^{[n+1]}=\{TI^{[n]}I^{[n]}\}+\{I^{[n]}TI^{[n]}\}+\{I^{[n]}I^{[n]}T\}.$$
\bprop
If $I$ is an ideal of a LeibTS $T$, then $I^{[n]}$ is also an ideal of $T$ and $I\supseteq I^{[1]}\supseteq\cdots\supseteq I^{[n]}\supseteq\cdots$, $\forall n\in\N$.
\eprop
\bpf
First we show $\{I^{[1]}TT\}\subseteq I^{[1]}$. By (\ref{VIP1}) and (\ref{VIP2}), it follows that
\begin{gather*}
\{\{TII\}TT\}\subseteq\!\{TI\{ITT\}\}\!+\!\{\{TIT\}IT\}\!+\!\{\{TIT\}TI\} \subseteq\!\{TII\}\!+\!\{IIT\}\!+\!\{ITI\}\subseteq\! I^{[1]},\\
\{\{ITI\}TT\}\subseteq\{IT\{ITT\}\}+\{\{ITT\}IT\}+\{\{ITT\}TI\} \subseteq\{ITI\}+\{IIT\}\subseteq I^{[1]},\\
\{\{IIT\}TT\}\subseteq\{\{ITI\}TT\}+\{\{ITT\}IT\}+\{I\{ITT\}T\} \subseteq I^{[1]}+\{IIT\}=I^{[1]},
\end{gather*}
which imply $\{I^{[1]}TT\}\subseteq I^{[1]}$. Similarly, one can prove $\{TI^{[1]}T\}\subseteq I^{[1]}$ and $\{TTI^{[1]}\}\subseteq I^{[1]}$. So $I^{[1]}$ is an ideal of $T$. Suppose that $I^{[k]}$ is an ideal of $T$. Then $I^{[k+1]}=(I^{[k]})^{[1]}$ is an ideal of $T$. Hence $I^{[n]}$ is an ideal of $T$ for all $n\in\N$. It is clear that $I\supseteq I^{[1]}\supseteq\cdots\supseteq I^{[n]}\supseteq\cdots$, $\forall n\in\N$.
\epf

\bdefn
An ideal $I$ of a LeibTS $T$ is called \textbf{solvable}, if there is a positive integer $k$ such that $I^{[k]}=0$.
\edefn

\bthm\label{sol eq}
Let $I$ be an ideal of a LeibTS $T$. Then the following statements are equivalent.
\benu[(1)]
\item $I$ is a solvable ideal of $T$.
\item $\I(I)$ is a solvable ideal of $U(T)$.
\item $\J(I)$ is a solvable ideal of $U(T)$.
\eenu
\ethm
\bpf
(1)$\Rightarrow$(2). We will use the induction to show
$$\I(I)^{(2n)}\subseteq I^{[n]}+I^{[n]}\otimes T.$$

Consider the base step $n=1$. We obtain
$$\I(I)^{(1)}=[\I(I),\I(I)]=[I+i(I),I+i(I)]=I\otimes I+[i(I),I]+[I,i(I)] +[i(I), i(I)].$$
Note that
\begin{gather*}
[i(I),I]=\{ITI\}+\{TII\}\subseteq I^{[1]},\quad
[I,i(I)]\subseteq \{IIT\}+\{ITI\}\subseteq I^{[1]},\\
\begin{aligned}
&[i(I), i(I)]=[I\otimes T+T\otimes I, I\otimes T+T\otimes I]\\
\subseteq&\{ITI\}\otimes T+\{ITT\}\otimes I+\{TII\}\otimes T+\{TIT\}\otimes I\subseteq I^{[1]}\otimes T+I\otimes I.
\end{aligned}
\end{gather*}
Hence $\I(I)^{(1)}\subseteq I^{[1]}+I^{[1]}\otimes T+I\otimes I$. So
\begin{align*}
&\I(I)^{(2)}=[\I(I)^{(1)}, \I(I)^{(1)}]\subseteq [I^{[1]}+I^{[1]}\otimes T+I\otimes I, I^{[1]}+I^{[1]}\otimes T+I\otimes I]\\
\subseteq&I^{[1]}\otimes I^{[1]}+\{I^{[1]}TI^{[1]}\}+\{III^{[1]}\} +\{I^{[1]}I^{[1]}T\}+\{I^{[1]}TI^{[1]}\}\otimes T+ \{I^{[1]}TT\}\otimes I^{[1]}\\
&+\{III^{[1]}\}\otimes T+\{IIT\}\otimes I^{[1]}+\{I^{[1]}II\}
+\{I^{[1]}TI\}\otimes I+\{III\}\otimes I\\
\subseteq&I^{[1]}+I^{[1]}\otimes T.
\end{align*}
Suppose $\I(I)^{(2k)}\subseteq I^{[k]}+I^{[k]}\otimes T$. Then
\begin{align*}
\I(I)^{(2k+1)}=&[\I(I)^{(2k)}, \I(I)^{(2k)}]\subseteq[I^{[k]} +I^{[k]}\otimes T, I^{[k]}+I^{[k]}\otimes T]\\
\subseteq&I^{[k]}\otimes I^{[k]}+\{I^{[k]}TI^{[k]}\}+\{I^{[k]}I^{[k]}T\} +\{I^{[k]}TI^{[k]}\}\otimes T+\{I^{[k]}TT\}\otimes I^{[k]}\\
\subseteq& I^{[k+1]}+I^{[k]}\otimes I^{[k]}+I^{[k+1]}\otimes T.
\end{align*}
Hence
\begin{align*}
&\I(I)^{(2k+2)}=[\I(I)^{(2k+1)}, \I(I)^{(2k+1)}]\\
\subseteq& [I^{[k+1]}+I^{[k]}\otimes I^{[k]}+I^{[k+1]}\otimes T, I^{[k+1]}+I^{[k]}\otimes I^{[k]}+I^{[k+1]}\otimes T]\\
\subseteq &I^{[k+1]}\otimes I^{[k+1]}+\{I^{[k]}I^{[k]}I^{[k+1]}\} +\{I^{[k+1]}TI^{[k+1]}\}+\{I^{[k+1]}I^{[k]}I^{[k]}\}+\!\{I^{[k]}I^{[k]}I^{[k]}\}\otimes I^{[k]}\\
&+\{I^{[k+1]}TI^{[k]}\}\otimes I^{[k]} +\{I^{[k+1]}I^{[k+1]}T\}+\!\{I^{[k]}I^{[k]}I^{[k+1]}\}\otimes T\\
&+\!\{I^{[k]}I^{[k]}T\}\otimes I^{[k+1]}+\!\{I^{[k+1]}TI^{[k+1]}\}\otimes T+\!\{I^{[k+1]}TT\}\otimes I^{[k+1]}\\
\subseteq&I^{[k+1]}+I^{[k+1]}\otimes T.
\end{align*}
Thus $\I(I)^{(2n)}\subseteq I^{[n]}+I^{[n]}\otimes T$ by induction, which implies that $\I(I)$ is solvable.

(2)$\Rightarrow$(1). It is sufficient to prove $I^{[n]}\subseteq \I(I)^{(n)}$ for all $n\in \N$. In fact,
\begin{align*}
I^{[1]}=&\{TII\}+\{ITI\}+\{IIT\}=[[T,I],I]+[[I,T],I]+[[I,I],T]\\
\subseteq&[[T,I],I]+[[I,T],I]+[I,[I,T]]\subseteq [i(I), I]+[I, i(I)] \subseteq[\I(I), \I(I)]=\I(I)^{(1)}.
\end{align*}
Suppose $I^{[k]}\subseteq \I(I)^{(k)}$. Then
\begin{align*}
I^{[k+1]}=&\{TI^{[k]}I^{[k]}\}+\{I^{[k]}TI^{[k]}\}+\{I^{[k]}I^{[k]}T\}\\
\subseteq&[[T,I^{[k]}],I^{[k]}]+[[I^{[k]},T],I^{[k]}]+[I^{[k]},[I^{[k]},T]]\\ \subseteq&[[T,\I(I)^{(k)}],\I(I)^{(k)}]+[[\I(I)^{(k)},T],\I(I)^{(k)}] +[\I(I)^{(k)},[\I(I)^{(k)},T]]\\
\subseteq&[\I(I)^{(k)},\I(I)^{(k)}]=\I(I)^{(k+1)}.
\end{align*}
Hence $I^{[n]}\subseteq \I(I)^{(n)}$ follows from the induction. So (1)$\Leftrightarrow$(2).

(3)$\Rightarrow$(2). It is clear since $\I(I)\subseteq\J(I)$.

(2)$\Rightarrow$(3). Note that
\begin{align*}
&[\J(I), U(T)]=[I+j(I), T+T\otimes T]\\
\subseteq&I\otimes T+[j(I), T]+\{ITT\}+[j(I), T]\otimes T\\
\subseteq&I\otimes T+I\subseteq\I(I).
\end{align*}
Then $\J(I)^{(1)}=[\J(I), \J(I)]\subseteq[\J(I), U(T)]\subseteq \I(I)$. Suppose $\J(I)^{(k)}\subseteq\I(I)^{(k-1)}$. Then
$$\J(I)^{(k+1)}=[\J(I)^{(k)}, \J(I)^{(k)}]\subseteq[\I(I)^{(k-1)}, \I(I)^{(k-1)}]=\I(I)^{(k)}.$$
It follows that $\J(I)^{(n+1)}\subseteq\I(I)^{(n)}$ by induction and so the solvability of $\I(I)$ implies the solvability of $\J(I)$. Thus (2)$\Leftrightarrow$(3).

Therefore, the above three items are equivalent and the proof is completed.
\epf

\bcor\label{sol equivalent}
Let $V=I+\M$ be a $\theta$-invariant ideal of $U(T)$. Then $I$ is a solvable ideal of $T$ if and only if $V$ is a solvable ideal of $U(T)$. In particular,
 $T$ is a solvable LeibTS if and only if $U(T)$ is a solvable Leibniz algebra.
\ecor

\bprop\label{sol ideal sum}
Suppose that $I$ and $J$ are solvable ideals of a LeibTS $T$. Then $I+J$ is a solvable ideal of $T$.
\eprop
\bpf
It is clear that $I+J$ is an ideal of $T$. Now we  use the induction to prove $(I+J)^{[n]}\subseteq I^{[n]}+J^{[n]}+I\cap J$. The base step $n=1$ is true, for
$$(I+J)^{[1]}=\{T,I+J,I+J\}+\{I+J,T,I+J\}+\{I+J,I+J,T\}\subseteq I^{[1]}+J^{[1]}+I\cap J.$$

Suppose $(I+J)^{[k]}\subseteq I^{[k]}+J^{[k]}+I\cap J$. Then
\begin{align*}
&(I+J)^{[k+1]}=\{T(I+J)^{[k]}(I+J)^{[k]}\} +\{(I+J)^{[k]}T(I+J)^{[k]}\} +\{(I+J)^{[k]}(I+J)^{[k]}T\}\\
\subseteq&\{T,I^{[k]}+J^{[k]}+I\cap J,I^{[k]}+J^{[k]}+I\cap J\} +\{I^{[k]}+J^{[k]}+I\cap J,T,I^{[k]}+J^{[k]}+I\cap J\} \\ &+\{I^{[k]}+J^{[k]}+I\cap J,I^{[k]}+J^{[k]}+I\cap J,T\}\\
\subseteq& I^{[k+1]}+J^{[k+1]}+I\cap J.
\end{align*}
So there is $n_1\in\N$ such that $(I+J)^{[n_1]}\subseteq I\cap J$. Note that $I\cap J$ is a solvable ideal of $T$, which implies $I+J$ is solvable.
\epf

\bdefn
Let $R(T)$ denote the sum of all solvable ideals in $T$. Then $R(T)$ is the unique maximal solvable ideal by Proposition \ref{sol ideal sum}, called the \textbf{solvable radical} of $T$.
\edefn

For a Leibniz algebra $L$, there is a corresponding Lie algebra $L/Ker(L)$, where $Ker(L)$ is a two-side ideal of $L$ spanned by $\{[x,x]\ |\ x\in L\}$. Similarly, for a LeibTS $T$, we define $Ker(T)$ to be the subspace of $T$ spanned by all elements in $T$ of the form $\{abc\}-\{acb\}+\{bca\}$, for all $a, b, c\in T$.

\bthm\label{Ker(T)}
Let $T$ be a LeibTS and $U(T)$ be its universal Leibniz envelope. Then
\benu[(1)]
\item $Ker(T)$ is an ideal of $T$ satisfying $\{TTKer(T)\}=\{TKer(T)T\}=0$.
\item $T/Ker(T)$ is a Lie triple system. Moreover, $T$ is a Lie triple system if and only if $Ker(T)=0$.
\item $Ker(T)\subseteq R(T)$.
\item $Ker(U(T))\cap T=Ker(T)$.
\item $\I(Ker(T))\subseteq Ker(U(T))\subseteq \J(Ker(T))$.
\eenu
\ethm
\bpf
(1) For $\{abc\}-\{acb\}+\{bca\}\in Ker(T)$ and $d, e\in T$, by (\ref{comm}) and (\ref{cyclic}), we have
\begin{gather*}
\{d,e,\{abc\}-\{acb\}+\{bca\}\}=\{d,e,\{abc\}+\{cab\}+\{bca\}\}=0,\\
\{d,\{abc\}-\{acb\}+\{bca\},e\}=\{d,\{abc\}+\{cab\}+\{bca\},e\}=0,
\end{gather*}
which imply $\{TTKer(T)\}=\{TKer(T)T\}=0$. Using (\ref{VIP1}) and (\ref{VIP2}), we have
\begin{align*}
\{\{abc\}-\{acb\},d,e\}=&\{\{abc\}de\}-\{\{acb\}de\}=\{a\{bcd\}e\}+\{\{adb\}ce\}-\{\{adc\}be\}\\
=&\{a\{bcd\}e\}+\{ad\{bce\}\}+\{\{ade\}bc\}-\{\{ade\}cb\},
\intertext{and}
\{\{bca\}de\}=&\{bc\{ade\}\}+\{\{bcd\}ae\}+\{\{bce\}ad\}-\{\{bce\}da\}.
\end{align*}
Note that
\begin{gather*}
\{\{ade\}bc\}-\{\{ade\}cb\}+\{bc\{ade\}\}\in Ker(T),\\
\{\{bce\}ad\}-\{\{bce\}da\}+\{ad\{bce\}\}\in Ker(T),
\end{gather*}
and
\begin{align*}
&\{a\{bcd\}e\}+\{\{bcd\}ae\}\\
=&(\{a\{bcd\}e\}-\{ae\{bcd\}\}+\{\{bcd\}ea\}) +(\{\{bcd\}ae\}-\{\{bcd\}ea\}+\{ae\{bcd\}\})\\
\in& Ker(T).
\end{align*}
Therefore, $\{\{abc\}-\{acb\}+\{bca\},d,e\}\in Ker(T)$ and $Ker(T)$ is an ideal of $T$.

(2) Since $Ker(T)$ is an ideal of $T$, $T/Ker(T)$ is also a LeibTS with the product $[\bar{a}\bar{b}\bar{c}]=\overline{\{abc\}}$.
$[\bar{a}\bar{a}\bar{b}]=\bar{0}$ comes from $\{aab\}=\{aab\}-\{aba\}+\{aba\}\in Ker(T)$, and then
$$[\bar{a}\bar{b}\bar{c}]+[\bar{b}\bar{c}\bar{a}]+[\bar{c}\bar{a}\bar{b}] =[\bar{a}\bar{b}\bar{c}]+[\bar{b}\bar{c}\bar{a}]-[\bar{a}\bar{c}\bar{b}] =\overline{\{abc\}-\{acb\}+\{bca\}}=\bar{0}.$$
Moreover,
\begin{align*}
[\bar{a}\bar{b}[\bar{c}\bar{d}\bar{e}]]=&\overline{\{ab\{cde\}\}} =\overline{\{\{abc\}de\}-\{\{abd\}ce\}-\{\{abe\}cd\}+\{\{abe\}dc\}}\\
=&\overline{\{\{abc\}de\}+\{c\{abd\}e\}+\{cd\{abe\}\}} =[[\bar{a}\bar{b}\bar{c}]\bar{d}\bar{e}]+[\bar{c}[\bar{a}\bar{b}\bar{d}]\bar{e}] +[\bar{c}\bar{d}[\bar{a}\bar{b}\bar{e}]].
\end{align*}
Hence $T/Ker(T)$ is a Lie triple system. It is clear $T$ is a Lie triple system if and only if $Ker(T)=0$.

(3) Since
$$Ker(T)^{[1]}=\{TKer(T)Ker(T)\}+\{Ker(T)TKer(T)\} +\{Ker(T)Ker(T)T\}=0,$$
$Ker(T)$ is a solvable ideal of $T$ and so $Ker(T)\subseteq R(T)$.

(4) $Ker(T)\subseteq Ker(U(T))\cap T$ follows from $\{abc\}-\{acb\}+\{bca\}=[a, b\otimes c]+[b\otimes c, a]$ $=[a+b\otimes c, a+b\otimes c]-[a,a]-[b\otimes c, b\otimes c]\in Ker(U(T))$.

Conversely, for all $x=[a+\sum b_i\otimes c_i, a+\sum b_i\otimes c_i]\in Ker(U(T))$, suppose $x=d+m$ with $d\in T, m\in T\otimes T$. Then $d=\sum (\{b_ic_ia\}+\{ab_ic_i\}-\{ac_ib_i\})\in Ker(T)$, that is, $Ker(U(T))\cap T\subseteq Ker(T)$.

(5) Note that $Ker(U(T))$ is an ideal of $U(T)$ and $Ker(U(T))\cap T=Ker(T)$. It follows from Theorem \ref{i(I)<I+M<j(I)} that $\I(Ker(T))\subseteq Ker(U(T))\subseteq \J(Ker(T))$.
\epf

Let $L$ be a Leibniz algebra and $Rad(L)$ the solvable radical of $L$. $L$ is called \textbf{semisimple} if $Rad(L)=Ker(L)$. Likewise, we introduce the definition of a semisimple LeibTS.

\bdefn
A LeibTS $T$ is said to be \textbf{semisimple}, if $R(T)=Ker(T)$.
\edefn

\bprop Let $T$ be a LeibTS. Then the following statements hold.
\benu[(1)]
\item $T/R(T)$ is a semisimple Lie triple system.
\item Let $I$ be an ideal of $T$. If $T/I$ is a semisimple LeibTS, then $R(T)\subseteq Ker(T)+I$.
\item If $I$ is an ideal of $T$ and $T/I$ is a semisimple Lie triple system, then $R(T)\subseteq I$.
\eenu
\eprop
\bpf
(1) $T/R(T)$ is clearly a Lie triple system since $Ker(T)\subseteq R(T)$. If $I/R(T)$ is a solvable ideal of $T/R(T)$, then $I^{[k]}/R(T)=(I/R(T))^{[k]}=\bar{0}$, for some $k\in\N$. Then $I^{[k]}\subseteq R(T)$, it follows that $I$ is a solvable ideal of $T$, so $I=R(T)$. Hence the only solvable ideal of $T/R(T)$ is $\bar{0}$ and $T/R(T)$ is semisimple.

(2) Since $T/I$ is semisimple, $R(T/I)=Ker(T/I)$. Note that $R(T)+I/I$ is solvable in $T/I$ and $Ker(T/I)=Ker(T)+I/I$. Then
$$ R(T)+I/I\subseteq R(T/I)=Ker(T/I)=Ker(T)+I/I,$$
which implies $R(T)\subseteq Ker(T)+I$.

(3) If $T/I$ is a Lie triple system, then $Ker(T/I)=0$, so $Ker(T)\subseteq I$, and hence $R(T)\subseteq I$ by using (2).
\epf

\bthm
$R(T)=Rad(U(T))\cap T$ and $Rad(U(T))=\I(R(T))+K(T)=\J(R(T))$.
\ethm
\bpf
Suppose that $\theta$ is the involutive automorphism of $U(T)$ determining $T$. Since $Rad(U(T))$ is solvable, its homomorphic image $\theta(Rad(U(T)))$ is solvable, then $\theta(Rad(U(T)))\subseteq Rad(U(T))$ and $Rad(U(T))$ is $\theta$-invariant. So
$$Rad(U(T))=Rad(U(T))\cap T+Rad(U(T))\cap(T\otimes T),$$ $$\I(Rad(U(T))\cap T)\subseteq Rad(U(T))\subseteq \J(Rad(U(T))\cap T).$$
Then $\I(Rad(U(T))\cap T)$ is solvable. By Theorem \ref{sol eq}, $Rad(U(T))\cap T$ and $\J(Rad(U(T))\cap T)$ is solvable, so $\J(Rad(U(T))\cap T)=Rad(U(T))$.

Since $Rad(U(T))\cap T$ is solvable, $Rad(U(T))\cap T\subseteq R(T)$, which implies $Rad(U(T))=\J(Rad(U(T))\cap T)\subseteq \J(R(T))$. But $\J(R(T))$ is solvable by the solvability of $R(T)$, then $\J(R(T))\subseteq Rad(U(T))$, hence $\J(R(T))=Rad(U(T))=\J(Rad(U(T))\cap T)$. So $R(T)=Rad(U(T))\cap T$ and $Rad(U(T))=\J(R(T))$.

Consider the factor Leibniz algebra $\overline{U(T)}=U(T)/\I(R(T))$. Since $\overline{U(T)}=\overline{T+T\otimes T}$ $=\overline{T}+\overline{T}\otimes \overline{T}=U(\overline{T})$ and $\overline{T}=T+\I(R(T))/\I(R(T))\cong T/(T\cap \I(R(T)))=T/R(T)$, it follows that $\overline{T}$ is a semisimple Lie triple system, so $R(\overline{T})=\overline{0}$. Hence
$$\overline{Rad(U(T))}\subseteq Rad(\overline{U(T)})=Rad(U(\overline{T})) =\J(R(\overline{T}))=\J(\overline{0})=\overline{\J(0)}=\overline{K(T)},$$
which implies $Rad(U(T))\subseteq K(T)+\I(R(T))$. On the other hand, note that $K(T)=\J(0)\!\subseteq\! \J(R(T))$ and $\I(R(T))\!\subseteq\! \J(R(T))$, so $Rad(U(T))\!=\!\I(R(T))\!+K(T)\!=\!\J(R(T))$.
\epf

\bthm
If $U(T)$ is a semisimple Leibniz algebra, then $T$ is a semisimple LeibTS; if $T$ is a semisimple LeibTS with $K(T)=0$, then $U(T)$ is a semisimple Leibniz algebra.
\ethm
\bpf
If $U(T)$ is semisimple, then $Rad(U(T))=Ker(U(T))$, so $$R(T)=Rad(U(T))\cap T=Ker(U(T))\cap T=Ker(T),$$
that is, $T$ is semisimple.

If $T$ is semisimple, then $R(T)=Ker(T)$, which implies $Rad(U(T))=\I(Ker(T))=\J(Ker(T))=Ker(U(T))$. Hence $U(T)$ is semisimple.
\epf

\section{Levi's theorem for Leibniz triple systems}
\bdefn{{\rm\cite{H}}}
A vector space $V$ is called a \textbf{module} for a Lie triple system $T$ if the vector space direct sum $T\dotplus V$ is itself a Lie triple system such that (1) $T$ is a subsystem, (2) $[abc]$ lies in $V$ if any one of $a, b, c$ is in $V$, (3) $[abc]=0$ if any two of $a, b, c$ are in $V$.
\edefn

\bdefn{{\rm\cite{Y}}}
Let $T$ be a Lie triple system and $V$ a vector space. $V$ is called a $T$-\textbf{module} if there exists a bilinear map $\delta: T\times T\rightarrow \End(V)$ such that for all $a,b,c,d\in T$,
\begin{gather*}
\delta(c, d)\delta(a,b)-\delta(b, d)\delta(a,c)-\delta(a, [bcd])+D(b, c)\delta(a,d)=0,\\
\delta(c, d)D(a, b)-D(a, b)\delta(c,d)+\delta([abc], d)+\delta(c, [abd])=0,
\end{gather*}
where $D(a,b)=\delta(b,a)-\delta(a,b)$.
\edefn

It is not difficult to prove that the above two definitions about the module for Lie triple systems are equivalent.

\bthm{{\rm\cite{H}}}\label{mod reduci}
If $T$ is a finite-dimensional semisimple Lie triple system over a field of characteristic 0, then every finite-dimensional module $V$ is completely reducible.
\ethm

\bthm[Levi's theorem for Lie triple systems]{{\rm\cite{L}}}
If $T$ is a Lie triple system, then there is a semisimple subsystem $S$ such that $T=R(T)\dotplus S$, where $R(T)$ is the solvable radical of $T$.
\ethm

\bthm[Levi's theorem for Leibniz triple systems]
Let $T$ be a LeibTS and $R(T)$ its solvable radical. Then there is a semisimple subsystem $S$ of $T$ such that $T=S\dotplus R(T)$. In particular, $S$ is a semisimple Lie triple system.
\ethm
\bpf
By Theorem \ref{Ker(T)}, $\overline{T}=T/Ker(T)$ is a Lie triple system. Then there are a semisimple subsystem $\overline{S}$ and the solvable radical $\overline{R}$ of $\overline{T}$ such that $\overline{T}=\overline{S}\dotplus \overline{R}$ from Levi's theorem for Lie triple systems. Let $\pi:T\rightarrow \overline{T}$ be the canonical projection. Then $\pi(R(T))$ is solvable since $R(T)$ is solvable, hence $\pi(R(T))\subseteq \overline{R}$.

Conversely, denote by $\pi^{-1}\left(\overline{R}\right)$ the inverse image of $\overline{R}$. Then $\pi^{-1}\left(\overline{R}\right)$ is a subsystem of $T$. Since $\Ker\pi=Ker(T)$, $\pi^{-1}\left(\overline{R}\right)/Ker(T)\cong \overline{R}$. Note that $Ker(T)$ and $\overline{R}$ are both solvable, which implies that $\pi^{-1}\left(\overline{R}\right)$ is solvable, then $\pi^{-1}\left(\overline{R}\right)\subseteq R(T)$, and so $\overline{R}\subseteq \pi(R(T))$. Thus $\pi(R(T))=\overline{R}$ and $R(T)=\pi^{-1}\left(\overline{R}\right)$.

Denote by $F=\pi^{-1}\left(\overline{S}\right)$ the inverse image of $\overline{S}$. Then $F$ is a subsystem of $T$ and $\overline{S}=F/Ker(T)$. Since $\overline{T}=T/Ker(T)$ and $\overline{R}=R(T)/Ker(T)$, it follows that
$$T/Ker(T)=F/Ker(T)\dotplus R(T)/Ker(T).$$
If $a\in F\cap R(T)$, then $\bar{a}\in (F/Ker(T))\cap (R(T)/Ker(T))=\overline{0}$, which implies $a\in Ker(T)$, hence $F\cap R(T)=Ker(T)$. So
$$T=F+R(T),\quad F\cap R(T)=Ker(T).$$

Define a bilinear map $\delta: \overline{S}\times \overline{S}\rightarrow \End(F)$ by $\delta(\bar{s_1}, \bar{s_2})(a)=\{as_1s_2\}$. Then $\{as_1s_2\}\in F$ since $\overline{S}=F/Ker(T)$ and $F$ is a subsystem of $T$. If $\bar{s_1}=\bar{s_1'}$ and $\bar{s_2}=\bar{s_2'}$, then $s_1-s_1', s_2-s_2'\in Ker(T)$, it follows from $\{TTKer(T)\}=\{TKer(T)T\}=0$ that
\begin{align*}
\{as_1s_2\}-\{as_1's_2'\}=&\{as_1s_2\}-\{as_1's_2\}+\{as_1's_2\}-\{as_1's_2'\}\\
=&\{a(s_1-s_1')s_2\}+\{as_1'(s_2-s_2')\}=0.
\end{align*}
Hence $\delta$ is well defined. Let $D(\bar{s_1}, \bar{s_2})(a)=\delta(\bar{s_2}, \bar{s_1})(a)-\delta(\bar{s_1}, \bar{s_2})(a)=\{as_2s_1\}-\{as_1s_2\}$. Then by (\ref{VIP2}) and (\ref{VIP3}),
\begin{align*}
&\delta(\bar{s_3}, \bar{s_4})\delta(\bar{s_1}, \bar{s_2})(a) -\delta(\bar{s_2}, \bar{s_4})\delta(\bar{s_1}, \bar{s_3})(a) -\delta(\bar{s_1}, \overline{\{s_2s_3s_4\}})(a) +D(\bar{s_2}, \bar{s_3})\delta(\bar{s_1}, \bar{s_4})(a)\\
=&\{\{as_1s_2\}s_3s_4\}-\{\{as_1s_3\}s_2s_4\}-\{as_1\{s_2s_3s_4\}\} +\{\{as_1s_4\}s_3s_2\}-\{\{as_1s_4\}s_2s_3\}\\
=&0,\\
&\delta(\bar{s_3}, \bar{s_4})D(\bar{s_1}, \bar{s_2})(a) -D(\bar{s_1}, \bar{s_2})\delta(\bar{s_3}, \bar{s_4})(a) +\delta(\overline{\{s_1s_2s_3\}}, s_4)(a) +\delta(\bar{s_3}, \overline{\{s_1s_2s_4\}})(a)\\
=&\{\{as_2s_1\}s_3s_4\}-\{\{as_1s_2\}s_3s_4\}-\{\{as_3s_4\}s_2s_1\} +\{\{as_3s_4\}s_1s_2\}\\
&+\{a\{s_1s_2s_3\}s_4\}+\{as_3\{s_1s_2s_4\}\}\\
=&0.
\end{align*}
Thus $F$ is an $\overline{S}$-module. By Theorem \ref{mod reduci}, $F$ is a completely reducible $\overline{S}$-module since $\overline{S}$ is a semisimple Lie triple system. Note that $Ker(T)$ is an ideal of $F$. Then $Ker(T)$ is a submodule of $F$. Hence by the completely reducibility of $F$, there is a complementary submodule $S$ of $F$ such that $F=S\dotplus Ker(T)$. Let $\sigma: F\rightarrow S$ be the canonical projection. Then $\Ker\sigma=Ker(T)$ and $S\cong F/Ker(T)$ is a semisimple Lie triple system such that
$$T=F+R(T)=S+Ker(T)+R(T)=S+R(T).$$

Similarly, $R(T)$ is an $\overline{S}$-module and there is a complementary submodule $R'$ of $R(T)$ such that $R(T)=R'\dotplus Ker(T)$. If $a\in S\cap R(T)$, then $a=a'+t$, for some $a'\in R'$ and $t\in Ker(T)$. Then $a'\in S\subseteq F$, which implies $a'\in R(T)\cap F=Ker(T)$. Moreover, $a'\in R'\cap Ker(T)=0$. Hence $a=t\in Ker(T)\cap S=0$, i.e., $R(T)\cap S=0$. Therefore, $T=S\dotplus R(T)$ with $S$ both a Lie triple system and a semisimple subsystem of $T$.
\epf

\section{The nilpotency of Leibniz triple systems}

For a Leibniz algebra $L$, the sequence of ideals
$$L\supseteq L^1\supseteq L^2\supseteq\cdots \quad \text{where~} L^{n+1}=[L^nL]$$
is called the central lower sequence of $L$.

Given an arbitrary LeibTS $T$, the central lower sequence of an ideal $I$ of $T$ is defined as
$I^{0}=I$ and for $n\geq1$,
$$I^{n+1}=\{I^{n}IT\}+\{I^{n}TI\}+\{II^{n}T\}+\{ITI^{n}\}+\{TII^{n}\}+\{TI^{n}I\}.$$
In particular, $T^0=T$ and $T^{n+1}=\{T^nTT\}+\{TT^nT\}+\{TTT^n\}$.

In \cite{AO}, the authors defined another sequence for a Leibniz algebra $L$ as $L^{<0>}=L$ and
$$L^{<n+1>}=[L^{<0>}L^{<n>}]+[L^{<1>}L^{<n-1>}]+\cdots+[L^{<n-1>}L^{<1>}]+[L^{<n>}L^{<0>}].$$
It is proved that $L^{<n>}=L^n$ for all $n\in\N$. Similarly, we could simplify the notation of $T^n$.

\bprop
For a LeibTS $T$, we define another sequence:
$$T^{(0)}=T,~ T^{(n+1)}=\{T^{(n)}TT\}.$$
Then $T^{(n)}=T^n$, for all $n\in\N$.
\eprop
\bpf
The proof is by induction on $n$. The base step $n=0$ is true, for $T^{(0)}=T^0=T$. Suppose $T^{(k)}=T^k$, $\forall k\leq n$. Note that
$$T^{n+1}=\{T^nTT\}+\{TT^nT\}+\{TTT^n\}=\{T^{(n)}TT\}+\{TT^{(n)}T\}+\{TTT^{(n)}\}.$$
So it is sufficient to prove
$$\{TT^{(n)}T\}+\{TTT^{(n)}\}\subseteq \{T^{(n)}TT\}=T^{(n+1)}.$$

In fact, the relation $\{T^{(i)}T^{(j)}T\}+\{T^{(i)}TT^{(j)}\}\subseteq T^{(i+j+1)}$ holds, for all $i, j\in\N$. We prove it by induction on $j$ for each $i$. If $j=0$, then
$$\{T^{(i)}T^{(0)}T\}+\{T^{(i)}TT^{(0)}\}=\{T^{(i)}TT\}=\{T^{(i+1)}TT\}, ~~\forall i\in\N.$$
Suppose $\{T^{(i)}T^{(j)}T\}+\{T^{(i)}TT^{(j)}\}\subseteq T^{(i+j+1)}, \forall i\in\N$. Then
\begin{align*}
 &\{T^{(i)}T^{(j+1)}T\}+\{T^{(i)}TT^{(j+1)}\}\\
=&\{T^{(i)}\{T^{(j)}TT\}T\}+\{T^{(i)}T\{T^{(j)}TT\}\}\\
\subseteq&\{\{T^{(i)}T^{(j)}T\}TT\}+\{\{T^{(i)}TT^{(j)}\}TT\}+\{\{T^{(i)}TT\}T^{(j)}T\}\\
&+\{\{T^{(i)}TT^{(j)}\}TT\}+\{\{T^{(i)}TT\}T^{(j)}T\}+\{\{T^{(i)}TT\}TT^{(j)}\}\\
\subseteq&\{T^{(i+j+1)}TT\}+\{T^{(i+1)}T^{(j)}T\}+\{T^{(i+1)}TT^{(j)}\}\\
=&T^{(i+j+2)}.
\end{align*}
Fix $i=0$, and it follows $\{TT^{(n)}T\}+\{TTT^{(n)}\}\subseteq T^{(n+1)}$, which completes the proof.
\epf

Hence the notions of $T^n$ and $T^{(n)}$ coincide.

\bdefn
An ideal $I$ of a LeibTS $T$ is called \textbf{nilpotent}, if there exists $n\in\N$ such that $I^n=0$.
\edefn

\bprop
If $I$ is an ideal of a LeibTS $T$, then $I^{n}$ is also an ideal of $T$ and $I^{n+1}\subseteq I^{n}$, $\forall n\in\N$.
\eprop
\bpf
We proceed by induction. It is clear when $n=0$. Suppose that $I^k$ is an ideal of $T$. Then
$$\{I^{k+1}TT\}=\{\{I^kTI\}+\{I^kIT\}+\{ITI^k\}+\{II^kT\}+\{TI^kI\}+\{TII^k\}, T, T\}.$$
By (\ref{VIP1}) and (\ref{VIP2}), it follows that
\begin{gather*}
\{\{I^kTI\}TT\}\subseteq\{I^kT\{ITT\}\}+\{\{I^kTT\}IT\}+\{\{I^kTT\}TI\} \subseteq I^{k+1},\\
\{\{I^kIT\}TT\}\subseteq\{I^k\{ITT\}T\}+\{\{I^kTI\}TT\}+\{\{I^kTT\}IT\} \subseteq I^{k+1},\\
\{\{ITI^k\}TT\}\subseteq \{IT\{I^kTT\}\}+\{\{ITT\}I^kT\}+\{\{ITT\}TI^k\} \subseteq I^{k+1},\\
\{\{II^kT\}TT\}\subseteq\{I\{I^kTT\}T\}+\{\{ITI^k\}TT\}+\{\{ITT\}I^kT\} \subseteq I^{k+1},\\
\{\{TI^kI\}TT\}\subseteq\{TI^k\{ITT\}\}+\{\{TI^kT\}IT\}+\{\{TI^kT\}TI\} \subseteq I^{k+1},\\
\{\{TII^k\}TT\}\subseteq\{TI\{I^kTT\}\}+\{\{TIT\}I^kT\}+\{\{TIT\}TI^k\} \subseteq I^{k+1},
\end{gather*}
which imply $\{I^{k+1}TT\}\subseteq I^{k+1}$. Similarly, one can prove
$$\{TI^{k+1}T\}\subseteq I^{k+1}\quad \text{and} \quad \{TTI^{k+1}\}\subseteq I^{k+1}.$$
So $I^{k+1}$ is an ideal of $T$. Hence $I^{n}$ is an ideal of $T$ and $I^{n+1}\subseteq I^{n}$, $\forall n\in\N$.
\epf

\bprop\label{nil ideal}
If $V$ is an ideal of $U(T)$, then $V\cap T$ is an ideal of $T$ such that $(V\cap T)^n\subseteq V^n\cap T$ for all $n\in\N$.
\eprop
\bpf
Since $V$ is an ideal of a Leibniz algebra $U(T)$, $V^0=V$ and $V^n=[V^{n-1}, V]$. Note that $V$ is an ideal of $U(T)$ if and only if $[V, U(T)]+[U(T), V]\subseteq V$. Then
\begin{gather*}
\{(V\cap T) T T\}=[[V\cap T, T], T]\subseteq[[V, U(T)], U(T)]\subseteq V,\\
\{T (V\cap T) T\}=[[T, V\cap T], T]\subseteq[[U(T), V], U(T)]\subseteq V,\\
\{T T (V\cap T)\}=[[T, T], V\cap T]\subseteq[U(T), V]\subseteq V,\\
\{(V\cap T) T T\}+\{T (V\cap T) T\}+\{T T (V\cap T)\}\subseteq\{TTT\}\subseteq T.
\end{gather*}
That is,
$$\{(V\cap T) T T\}+\{T (V\cap T) T\}+\{T T (V\cap T)\}\subseteq V\cap T.$$
So $V\cap T$ is an ideal of $T$.

Now we use the induction to prove $(V\cap T)^n\subseteq V^n\cap T$ for all $n\in\N$. The base step $n=0$ is clearly true. Suppose $(V\cap T)^k\subseteq V^k\cap T$. Then
\begin{align*}
&(V\cap T)^{k+1}\\
=&\{(V\cap T)^k(V\cap T)T\}+\{(V\cap T)^kT(V\cap T)\}+\{(V\cap T)(V\cap T)^kT\}\\
 &+\{(V\cap T)T(V\cap T)^k\}+\{T(V\cap T)(V\cap T)^k\}+\{T(V\cap T)^k(V\cap T)\}\\
\subseteq&\{(V^k\cap T)(V\cap T)T\}+\{(V^k\cap T)T(V\cap T)\}+\{(V\cap T)(V^k\cap T)T\}\\
 &+\{(V\cap T)T(V^k\cap T)\}+\{T(V\cap T)(V^k\cap T)\}+\{T(V^k\cap T)(V\cap T)\}\\
\subseteq&[[V^k, V], T]+[[V^k, T], V]+[[V, V^k], T]+[[V, T], V^k]+[[T, V], V^k]+[[T, V^k], V]\\
\subseteq&[[V^k, T], V]+[V^k, [V, T]]+[[V, T], V^k]+[V, [V^k, T]]+[V, V^k]+[V^k, V]\\
\subseteq&[V^k, V]+[V, V^k]=V^{k+1}.
\end{align*}
Then $(V\cap T)^{k+1}\subseteq V^{k+1}\cap T$. Therefore, $(V\cap T)^n\subseteq V^n\cap T$ for all $n\in\N$.
\epf

\blem\label{I&I(I)}
Let $I$ be an ideal of a LeibTS $T$. Then

(1) $I^n\subseteq \I(I)^n$, for all $n\in\N$.

(2) $\I(I)^{2n-1}\subseteq I^{n}+I^{n-1}\otimes I+I^{n}\otimes T$ and $\I(I)^{2n}\subseteq I^{n}+I^{n}\otimes T$, $\forall n\geq1$.
\elem
\bpf
(1) It follows from Proposition \ref{nil ideal}.

(2) We still use the induction. It is straightforward that
\begin{align*}
\I(I)^1=&[\I(I), \I(I)]=[I+I\otimes T+T\otimes I, I+I\otimes T+T\otimes I]\\
\subseteq&I\!\otimes\! I\!+\!\{IIT\}\!+\!\{ITI\}\!+\!\{ITI\}\!\otimes\! T\!+\!\{ITT\}\!\otimes\! I\!+\!\{TII\}\!+\!\{TII\}\!\otimes\! T\!+\!\{TIT\}\!\otimes\! I\\
\subseteq& I^1+I\otimes I+I^1\otimes T
\intertext{and}
\I(I)^2=&[\I(I)^1, \I(I)]\subseteq[I^1+I\otimes I+I^1\otimes T, I+I\otimes T+T\otimes I]\\
\subseteq&I^1\!\otimes\! I\!+\!\{I^1IT\}\!+\!\{I^1TI\}\!+\!\{III\}\!+\!\{III\}\!\otimes\! T\!+\!\{IIT\}\!\otimes\! I\!+\!\{I^1TI\}\!\otimes\! T\!+\!\{I^1TT\}\!\otimes\! I\\
\subseteq&I^1+I^1\otimes T.
\end{align*}

Suppose $\I(I)^{2k-1}\subseteq I^{k}+I^{k-1}\otimes I+I^{k}\otimes T$ and $\I(I)^{2k}\subseteq I^{k}+I^{k}\otimes T$. Then
\begin{align*}
\I(I)^{2k+1}=&[\I(I)^{2k}, \I(I)]\subseteq[I^{k}+I^{k}\otimes T, I+I\otimes T+T\otimes I]\\
\subseteq&I^k\!\otimes\! I\!+\!\{I^kIT\}\!+\!\{I^kTI\}\!+\!\{I^kTI\}\!\otimes\! T\!+\!\{I^kTT\}\!\otimes\! I\\
\subseteq& I^{k+1}+I^k\otimes I+I^{k+1}\otimes T
\intertext{and}
\I(I)^{2k+2}=&[\I(I)^{2k+1}, \I(I)]\subseteq[I^{k+1}+I^k\otimes I+I^{k+1}\otimes T, I+I\otimes T+T\otimes I]\\
\subseteq&I^{k+1}\otimes I+\{I^{k+1}IT\}+\{I^{k+1}TI\}+\{I^kII\}+\{I^kII\}\otimes T+\{I^kIT\}\otimes I\\
        &+\{I^{k+1}TI\}\otimes T+\{I^{k+1}TT\}\otimes I\\
\subseteq&I^{k+1}+I^{k+1}\otimes T.
\end{align*}
The proof is completed.
\epf

\blem\label{nil ideal sum}
Let $I$ and $J$ be two ideals of a LeibTS $T$. Then
$$(I+J)^n\subseteq I^n+J^n+\sum_{k=0}^{n-1}I^k\cap J^{n-1-k},\quad \forall n\geq1.$$
\elem
\bpf
If $n=1$, then
$$(I+J)^1=\{(I+J)(I+J)T\}+\{(I+J)T(I+J)\}+\{T(I+J)(I+J)\}\subseteq I^1+J^1+I\cap J.$$
Suppose $(I+J)^m\subseteq I^m+J^m+\sum_{k=0}^{m-1}I^k\cap J^{m-1-k}$. Then
\begin{align*}
(I+J)^{m+1}=&\{(I+J)^m(I+J)T\}+\{(I+J)^mT(I+J)\}+\{(I+J)(I+J)^mT\}\\
            &+\{(I+J)T(I+J)^m\}+\{T(I+J)(I+J)^m\}+\{T(I+J)^m(I+J)\}
\end{align*}
Note that
\begin{align*}
&\{(I+J)^m(I+J)T\}\subseteq \{I^m+J^m+\sum_{k=0}^{m-1}I^k\cap J^{m-1-k}, I+J, T\}\\
                 \subseteq& I^{m+1}+I^m\cap J+I\cap J^m+J^{m+1}+\sum_{k=0}^{m-1}I^{k+1}\cap J^{m-1-k}+\sum_{k=0}^{m-1}I^k\cap J^{m-k}\\
                 =&I^{m+1}+J^{m+1}+\sum_{k=0}^{m}I^k\cap J^{m-k}.
\end{align*}

One can prove that other five items belong to  $I^{m+1}+J^{m+1}+\sum_{k=0}^{m}I^k\cap J^{m-k}$ in the same way. Hence the lemma follows from the induction.
\epf

By Lemmas \ref{I&I(I)} and \ref{nil ideal sum}, it is clear that the following theorem holds.
\bthm\label{I and I(I)}
Let $T$ be a LeibTS. Then

(1) $I$ is a nilpotent ideal of $T$ if and only if $\I(I)$ is a nilpotent ideal of $U(T)$. In particular, $T$ is a nilpotent LeibTS if and only if $U(T)$ is a nilpotent Leibniz algebra.

(2) If $I$ and $J$ are nilpotent ideals of $T$, then $I+J$ is a nilpotent ideal of $T$.
\ethm

\bdefn
Let $N(T)$ denote the sum of all nilpotent ideals in $T$. Then $N(T)$ is the unique maximal nilpotent ideal by Theorem \ref{I and I(I)}, called the \textbf{nilpotent radical} of $T$.
\edefn

\bthm
Suppose $Nil(U(T))$ is the nilpotent radical of $U(T)$. Then $N(T)=Nil(U(T))\cap T$.
\ethm
\bpf
Since $Nil(U(T))$ is a nilpotent ideal of $U(T)$, $Nil(U(T))\cap T$ is a nilpotent ideal of $T$ by Proposition \ref{nil ideal}. If $I$ is a nilpotent ideal of $T$, then $\I(I)$ a nilpotent ideal of $U(T)$ by Theorem \ref{I and I(I)}, hence $\I(I)\subseteq Nil(U(T))$. So $I=\I(I)\cap T\subseteq Nil(U(T))\cap T$. Therefore, $Nil(U(T))\cap T$ is a maximal nilpotent ideal of $T$, that is, $Nil(U(T))\cap T=N(T)$.
\epf

\blem{{\rm\cite{AO}}}
If $L$ is a solvable Leibniz algebra, then $[L, L]$ is nilpotent.
\elem

\bprop
Suppose $T$ is a LeibTS. If $T$ is nilpotent, then $T$ is solvable; if $T$ is solvable, then $T^1=\{TTT\}$ is nilpotent.
\eprop
\bpf
It is clear that $T^{[n]}\subseteq T^n$. So the nilpotency of $T$ implies the solvability of $T$. Conversely, if $T$ is solvable, then $U(T)$ is solvable by Corollary \ref{sol equivalent}, so $[U(T), U(T)]$ is nilpotent, hence $[U(T), U(T)]\subseteq Nil(U(T))$. Note that
$$T^1=\{TTT\}=[[T, T], T]\subseteq [U(T), U(T)]\subseteq Nil(U(T)).$$
Then $T^1\subseteq Nil(U(T))\cap T=N(T)$. So $T^1$ is nilpotent.
\epf

\section{The representation of Leibniz triple systems}

\bdefn
Let $T$ be a LeibTS and $V$ a vector space. $V$ is called a $T$-\textbf{module}, if $T\dotplus V$ is a LeibTS such that (1) $T$ is a subsystem, (2) $\{abc\}\in V$ if any one of $a, b, c$ is in $V$, (3) $\{abc\}=0$ if any two of $a, b, c$ are in $V$.
\edefn

\bdefn
Let $T$ be a LeibTS and $V$ a vector space. Suppose $l, m, r: T\times T\rightarrow \End(V)$ are bilinear maps such that
\begin{align}
l(a, \{bcd\})&=l(\{abc\}, d)-l(\{acb\}, d)-l(\{adb\}, c)+l(\{adc\}, b), \label{represen1}\\
m(a, d)l(b, c)&=m(\{abc\}, d)-m(\{acb\}, d)-r(c, d)m(a, b)+r(b, d)m(a, c), \label{represen2}\\
m(a, d)m(b, c)&=r(c, d)l(a, b)-r(c, d)m(a, b)-m(\{acb\}, d)+r(b, d)l(a, c), \label{represen3}\\
m(a, d)r(b, c)&=r(c, d)m(a, b)-r(c, d)l(a, b)-r(b, d)l(a, c)+m(\{acb\}, d), \label{represen4}\\
r(\{abc\}, d)&=r(c, d)r(a, b)-r(c, d)r(b, a)-r(b, d)r(c, a)+r(a, d)r(c, b), \label{represen5}\\
l(a, b)l(c, d)&=l(\{abc\}, d)-l(\{abd\}, c)-r(c, d)l(a, b)+r(d, c)l(a, b), \label{represen6}\\
l(a, b)m(c, d)&=m(\{abc\}, d)-r(c, d)l(a, b)-l(\{abd\}, c)+m(\{abd\}, c), \label{represen7}\\
l(a, b)r(c, d)&=r(c, d)l(a, b)-m(\{abc\}, d)-m(\{abd\}, c)+l(\{abd\}, c), \label{represen8}\\
m(a, \{bcd\})&=r(c, d)m(a, b)-r(b, d)m(a, c)-r(b, c)m(a, d)+r(c, b)m(a, d), \label{represen9}\\
r(a, \{bcd\})&=r(c, d)r(a, b)-r(b, d)r(a, c)-r(b, c)r(a, d)+r(c, b)r(a, d). \label{represen10}
\end{align}
Then $(r, m, l)$ is called a \textbf{representation} of $T$ on $V$.
\edefn

\bre
If $V$ is a $T$-module and define $l(a, b)(v)=\{abv\}$, $m(a, b)(v)=\{avb\}$, $r(a, b)(v)=\{vab\}$, then $(r, m, l)$ is a representation of $T$. Conversely, if $(r, m, l)$ is a representation of $T$, then $V$ is a $T$-module by setting $\{abv\}=l(a, b)(v)$, $\{avb\}=m(a, b)(v)$ and $\{vab\}=r(a, b)(v)$. So the notions of a $T$-module and a representation of $T$ are equivalent.
\ere

\bprop
Suppose that $T$ is a LeibTS and $(r, m, l)$ is a representation of $T$ on $V$.

(1) $m(a, b)^{k+1}=(-1)^km(a, b)r(a, b)^k, \forall k\in\N$. So if $r(a, b)$ is nilpotent, then $m(a, b)$ is nilpotent.

(2) Set $R(a, b)=r(a, b)-r(b, a)$. Then $l(a, b)^{k+1}=(-1)^kl(a, b)R(a, b)^k, \forall k\in\N$. So if $R(a, b)$ is nilpotent, then $l(a, b)$ is nilpotent.
\eprop
\bpf
(1) By (\ref{represen3}) and (\ref{represen4}), $m(a, b)m(c, d)=-m(a, b)r(c, d)$. Then
$$m(a, b)^{k+1}=-m(a, b)^{k}r(a, b)=\cdots=(-1)^km(a, b)r(a, b)^k, ~~\forall k\in\N.$$

(2) By (\ref{represen8}),
\begin{gather*}
l(a, b)r(c, d)=r(c, d)l(a, b)-m(\{abc\},  d)-m(\{abd\},  c)+l(\{abd\},  c),\\
l(a, b)r(d, c)=r(d, c)l(a, b)-m(\{abd\},  c)-m(\{abc\},  d)+l(\{abc\},  d).
\end{gather*}
Then
$$l(a, b)R(c, d)=-l(\{abc\}, d)+l(\{abd\}, c)+r(c, d)l(a, b)-r(d, c)l(a, b)=-l(a, b)l(c, d),$$
where the last identity is by (\ref{represen6}). Hence
$$l(a, b)^{k+1}=-l(a, b)^{k}R(a, b)=\cdots=(-1)^kl(a, b)R(a, b)^k, ~~\forall k\in\N.$$
The proof is completed.
\epf

Denote by $l(T, T)=\{\sum l(a, b)\ |\ a, b\in T\}$, $m(T, T)=\{\sum m(a, b)\ |\ a, b\in T\}$, $r(T, T)=\{\sum r(a, b)\ |\ a, b\in T\}$ and $R(T, T)=\{\sum R(a, b)\ |\ a, b\in T\}$. Then $l(T, T)$, $m(T, T)$, $r(T, T)$ and $R(T, T)$ are subspaces of $\End(V)$. Note that for all $a, b, c, d\in T$, by (\ref{represen5}) and (\ref{represen10}),
\begin{align*}
 &[R(a, b), R(c, d)]=R(a, b)R(c, d)-R(c, d)R(a, b)\\
=&(r(a, b)-r(b, a))(r(c, d)-r(d, c))-(r(c, d)-r(d, c))(r(a, b)-r(b, a))\\
=&r(a, b)r(c, d)-r(b, a)r(c, d)-r(a, b)r(d, c)+r(b, a)r(d, c)\\
 &-r(c, d)r(a, b)+r(d, c)r(a, b)+r(c, d)r(b, a)-r(d, c)r(b, a)\\
=&r(c, \{dab\})+r(d, b)r(c, a)+r(d, a)r(c, b)-r(a, d)r(c, b)\\
 &-r(c, \{dba\})-r(d, a)r(c, b)-r(d, b)r(c, a)+r(b, d)r(c, a)\\
 &-r(d, \{cab\})-r(c, b)r(d, a)-r(c, a)r(d, b)+r(a, c)r(d, b)\\
 &+r(d, \{cba\})+r(c, a)r(d, b)+r(c, b)r(d, a)-r(b, c)r(d, a)\\
 &-r(\{cba\}, d)+r(a, d)r(c, b)-r(a, d)r(b, c)-r(b, d)r(a, c)\\
 &+r(\{dba\}, c)-r(a, c)r(d, b)+r(a, c)r(b, d)+r(b, c)r(a, d)\\
 &+r(\{cab\}, d)-r(b, d)r(c, a)+r(b, d)r(a, c)+r(a, d)r(b, c)\\
 &-r(\{dab\}, c)+r(b, c)r(d, a)-r(b, c)r(a, d)-r(a, c)r(b, d)\\
=&R(c, \{dab\})-R(c, \{dba\})-R(d, \{cab\})+R(d, \{cba\}).
\end{align*}
So $R(T, T)$ is a subalgebra of $gl(V)$.

\bprop
$l(T, T)$, $m(T, T)$ and $r(T, T)$ are all adjoint $R(T, T)$-module.
\eprop
\bpf
It is sufficient to prove that
$[R(T, T), l(T, T)]\subseteq l(T, T)$, $[R(T, T), m(T, T)]\subseteq m(T, T)$ and $[R(T, T), r(T, T)]\subseteq r(T, T)$.

It follows from (\ref{represen8}) that
\begin{gather*}
[r(c, d), l(a, b)]=m(\{abc\}, d)+m(\{abd\}, c)-l(\{abd\}, c),\\
[r(d, c), l(a, b)]=m(\{abd\}, c)+m(\{abc\}, d)-l(\{abc\}, d).
\end{gather*}
Then $[R(c, d), l(a, b)]=l(\{abc\}, d)-l(\{abd\}, c)\in l(T, T)$.

It follows from (\ref{represen4}) that
\begin{gather*}
m(a, d)r(b, c)=r(c, d)m(a, b)-r(c, d)l(a, b)-r(b, d)l(a, c)+m(\{acb\}, d),\\
m(a, d)r(c, b)=r(b, d)m(a, c)-r(b, d)l(a, c)-r(c, d)l(a, b)+m(\{abc\}, d).
\end{gather*}
Then $m(a, d)R(b, c)=r(c, d)m(a, b)-r(b, d)m(a, c)+m(\{acb\}, d)-m(\{abc\}, d)$. Note that by (\ref{represen9}),
$$R(b, c)m(a, d)=r(c, d)m(a, b)-r(b, d)m(a, c)-m(a, \{bcd\}).$$
Hence $[R(b, c), m(a, d)]=-m(a, \{bcd\})-m(\{acb\}, d)+m(\{abc\}, d)\in m(T, T)$.

It follows from (\ref{represen5}) and (\ref{represen10}) that
\begin{gather*}
r(c, d)R(a, b)=r(\{abc\}, d)+r(b, d)r(c, a)-r(a, d)r(c, b),\\
R(a, b)r(c, d)=-r(c, \{abd\})+r(b, d)r(c, a)-r(a, d)r(c, b).
\end{gather*}
Then $[R(a, b), r(c, d)]=-r(c, \{abd\})-r(\{abc\}, d)\in r(T, T)$.
\epf

Set
$$Z_T(V)=\{x\in T\ |\ \{xTV\}=\{xVT\}=\{TxV\}=\{VxT\}=\{VTx\}=\{TVx\}=0\}.$$
If $T$ is a Lie triple system, then it is clear that $Z_T(V)$ is an ideal of $T$. But it is not true if $T$ is a non-Lie LeibTS.

In fact, it is easy to prove that
\begin{gather*}
\{TV\{Z_T(V)TT\}\}=\{TV\{TZ_T(V)T\}\}=\{TV\{TTZ_T(V)\}\}=0,\\
\{VT\{Z_T(V)TT\}\}=\{VT\{TZ_T(V)T\}\}=\{VT\{TTZ_T(V)\}\}=0,\\
\{T\{Z_T(V)TT\}V\}=\{T\{TZ_T(V)T\}V\}=\{T\{TTZ_T(V)\}V\}=0,\\
\{V\{Z_T(V)TT\}T\}=\{V\{TZ_T(V)T\}T\}=\{V\{TTZ_T(V)\}T\}=0,\\
\{\{TTZ_T(V)\}TV\}=\{\{TZ_T(V)T\}TV\}=\{\{TTZ_T(V)\}VT\}=\{\{TZ_T(V)T\}VT\}=0.
\end{gather*}
However, we fail to show $\{\{Z_T(V)TT\}TV\}=0$ or $\{\{Z_T(V)TT\}VT\}=0$. Therefore, one can not say that $Z_T(V)$ is an ideal of $T$.
\bthm
Suppose that $T$ is a LeibTS and $V$ is an irreducible $T$-module.
If $Z_T(V)$ is an ideal of $T$, then $T/Z_T(V)$ is a Lie triple system, and either $l(T, T)=m(T, T)=0$ or
$$l(x, y)=-l(y, x)=R(x, y),\quad r(x, y)=-m(x, y),\quad \forall x, y\in T.$$
\ethm
\bpf
Let $\tilde{T}$ be a split extension of $V$ by $T$. Then $\tilde{T}=T\dotplus V$ and $V$ is also a $\tilde{T}$-module. Note that
\begin{align*}
Z_{\tilde{T}}(V)&=\{x\in \tilde{T}\ |\ \{x\tilde{T}V\}=\{xV\tilde{T}\}=\{\tilde{T}xV\}=\{Vx\tilde{T}\}=\{V\tilde{T}x\}=\{\tilde{T}Vx\}=0\}\\
&=\{x\in \tilde{T}\ |\ \{xTV\}=\{xVT\}=\{TxV\}=\{VxT\}=\{VTx\}=\{TVx\}=0\}\\
&=Z_T(V)\dotplus V.
\end{align*}

Consider the following two cases.

1. $Z_T(V)=0$. Since $V$ is an irreducible $T$-module, $V$ is a minimal ideal of $\tilde{T}$. Note that $Ker(\tilde{T})$ is an ideal of $\tilde{T}$. Then either $V\subseteq Ker(\tilde{T})$ or $V\cap Ker(\tilde{T})=0$. Then
$\{V\tilde{T}Ker(\tilde{T})\}=\{VKer(\tilde{T})\tilde{T}\}=\{\tilde{T}VKer(\tilde{T})\}=\{Ker(\tilde{T})V\tilde{T}\}=\{\tilde{T}Ker(\tilde{T})V\}=\{Ker(\tilde{T})\tilde{T}V\}=0$ by Theorem \ref{Ker(T)}.
Hence $Ker(\tilde{T})\subseteq Z_{\tilde{T}}(V)=V$, and so $Ker(T)\subseteq Ker(\tilde{T})\subseteq V$. Then $Ker(T)\subseteq V\cap T=0$ and $\tilde{T}/V\cong T$ is a Lie triple system. It follows from $Ker(\tilde{T})\subseteq V$ and the minimality of $V$ that $Ker(\tilde{T})=0$ or $Ker(\tilde{T})=V$.

If $Ker(\tilde{T})=V$, then $\{\tilde{T}\tilde{T}V\}=\{\tilde{T}V\tilde{T}\}=0$, which implies $\{TTV\}=\{TVT\}=0$, that is, $l(T, T)=m(T, T)=0$.

If $Ker(\tilde{T})=0$, then $\tilde{T}$ is a Lie triple system. Hence $l(x, y)=-l(y, x)=R(x, y)$ and $r(x, y)=-m(x, y)$, for all $x, y\in T$.

2. $Z_T(V)\neq0$. Let $\bar{T}=T/Z_T(V)$. Then $\bar{T}$ is a LeibTS and $V$ is a $\bar{T}$-module satisfying $Z_{\bar{T}}(V)=0$. Note that $V$ is also an irreducible $\bar{T}$-module. Hence $\bar{T}$ is a Lie triple system.
\epf

\bthm\label{char vector}
Let $T$ be a LeibTS and $V\neq0$ a finite-dimensional $T$-module.  Suppose that $Z_T(V_0)$ is an ideal of $T$ for every submodule $V_0$ of $V$. If $R(x, y)$ is nilpotent for all $x, y\in T$, then $l(x, y)$ is nilpotent and there exists $0\neq v\in V$ such that
$$R(x, y)(v)=l(x, y)(v)=0,\quad \forall x, y\in T.$$
\ethm
\bpf
Since $V$ is finite-dimensional, there exists a sequence of submodules
$$V=V_0\supset V_1\supset\cdots\supset V_{n-1}\supset V_n=0$$
such that each factor $V_k/V_{k+1}$ ($k=0,1, \cdots, n-1$) is an irreducible $T$-module. Then
$$V\cong V_0/V_1\dotplus V_1/V_2\dotplus \cdots\dotplus V_{n-2}/V_{n-1}\dotplus V_{n-1}/V_n.$$

Note that $R(T, T)$ is a subalgebra of $gl(V)$ and by the Engel's theorem for Lie algebras, there exists $0\neq v\in V$ such that $R(x, y)(v)=0$, $\forall x, y\in T$. Suppose $v=\overline{v_1}+\cdots+\overline{v_{n-1}}$, where $\overline{v_i}\in V_i/V_{i+1}$. Since $V_i/V_{i+1}$ is an irreducible $T$-module, it follows that $l(x, y)|_{V_i/V_{i+1}}=0$ or $l(x, y)|_{V_i/V_{i+1}}=-R(y, x)|_{V_i/V_{i+1}}$. In either case we have $l(x, y)(\overline{v_i})=0$, so $l(x, y)(v)=0$, $\forall x, y\in T$.
\epf

\bre
If $Z_T(I)$ is an ideal of $T$ for every ideal $I$ of $T$, then by Theorem \ref{char vector}, there exits a basis of a finite-dimensional LeibTS $T$ relative to which the matrix of $l(x, y)$ is strictly upper triangular for all $x, y\in T$. Then $\{TTT^n\}=0$ for some $n\in\N$. But one can not say that $T$ is nilpotent. Hence the Engel's theorem for a LeibTS could not be proved in this way.
\ere


\begin{thebibliography}{99}
\bibitem{AAO} S. Albeverio, S. Ayupov, B. Omirov, (2005), On nilpotent and simple Leibniz algebras. Comm. Algebra 33(1), 159-172.
\bibitem{AAO2} S. Albeverio, S. Ayupov, B. Omirov, (2006), Cartan subalgebras, weight spaces, and criterion of solvability of finite dimensional Leibniz algebras. Rev. Mat. Complut. 19(1), 183-195.
\bibitem{AO} S. Ayupov, B. Omirov, (1998), On Leibniz algebras. Algebra and operator theory, Kluwer Acad. Publ., Dordrecht, 1-12.
\bibitem{B3} D. Barnes, (2012), On Engel's theorem for Leibniz algebras. Comm. Algebra 40(4), 1388-1389.
\bibitem{B4} D. Barnes, (2012), On Levi's theorem for Leibniz algebras. Bull. Aust. Math. Soc. 86(2), 184-185.
\bibitem{B2} A. Bolh, (1965), On a generalization of the concept of Lie algebra. (Russian) Dokl. Akad. Nauk SSSR 165, 471-473.
\bibitem{BS} M. Bremner, J. S\'{a}nchez-Ortega, (2014), Leibniz triple systems. Commun. Contemp. Math. 16(1), 1350051, 19 pp.
\bibitem{DMS} I. Demir, K. Misra, E. Stitzinger, (2014), On some structures of Leibniz algebras. Contemp. Math. 623, Amer. Math. Soc., Providence, RI.
\bibitem{GKO} S. G\'{o}mez-Vidal, S. Khudoyberdiyev, B. Omirov, (2014), Some remarks on semisimple Leibniz algebras. J. Algebra 410, 526-540.
\bibitem{G} V. Gorbatsevich, (2013), On some basic properties of Leibniz algebras. arXiv: 1302. 3345.
\bibitem{H} B. Harris, (1961), Cohomology of Lie triple systems and Lie algebras with involution. Trans. Amer. Math. Soc. 98, 148-162.
\bibitem{K} P. Kolesnikov, (2008), Varieties of dialgebras, and conformal algebras. (Russian) Sibirsk. Mat. Zh. 49(2), 322-339; translation in Sib. Math. J. 49(2), 257-272.
\bibitem{L} W. Lister, (1952), A structure theory of Lie triple systems. Trans. Amer. Math. Soc. 72, 217-242.
\bibitem{L2} J. Loday, (1993), Une version non commutative des alg\`{e}bres de Lie: les alg\`{e}bres de Leibniz. (French) Enseign. Math. 39(3-4), 269-293.
\bibitem{LP} J. Loday, T. Pirashvili, (1993), Universal enveloping algebras of Leibniz algebras and (co)homology. Math. Ann. 296(1), 139-158.
\bibitem{M} K. Meyberg, (1972), Lectures on algebras and triple systems. The University of Virginia, Charlottesville, Va. Notes on a course of lectures given during the academic year 1971-1972.
\bibitem{P} A. Patsourakos, (2007), On nilpotent properties of Leibniz algebras. Comm. Algebra 35(12), 3828-3834.
\bibitem{Y} K. Yamaguti, (1960), On the cohomology space of Lie triple system. Kumamoto J. Sci. Ser. A 5, 44-52.
\end{thebibliography}
\end{document}